\documentclass[12pt]{article}
\usepackage{amsmath}
\usepackage{amssymb}
\usepackage{amsfonts}
\usepackage{mitpress}

\newdimen\dummy
\dummy=\oddsidemargin
\input{tcilatex}
\begin{document}

\title{Frierson's 1907 Parameterization of Compound Magic Squares Extended
to Orders $3^{l}$, $l=1,2,3,..$, with Information Entropy}
\author{P.D.Loly and I.D.Cameron \\
Department of Physics and Astronomy, University of Manitoba, Winnipeg,
Manitoba, R3T 2N2, Canada}
\maketitle

\begin{abstract}
Frierson used a powerful parameterization of the pattern of the order $3$
associative magic square to construct a family of six related order \textbf{%
\ }$3^{2}=9$ compound (or composite) magic squares, several of them ancient.
Stimulated by Bellew's 1997 extension to order $27=3^{3}$, we extend these
ideas to all orders $3^{l}$, $l=1,2,3,..$, and in addition find simple
formulae for the matrix spectra and entropic measures for all those orders.
This construction is fractal and we give numerical results to order $%
243=3^{6}$ which show an information entropy measure converging to a
constant value of about $1.168.$. for the lowest entropy members.

We also briefly consider compounding of an order $4$ magic square with the
lowest entropy, for which we find a similar trend to constant entropy.
\end{abstract}

\section{Introduction\label{intro}}

Magic squares (MSs) have the same line sum for all Rows, Columns, and their
two main Diagonals (RCDs), with most interest in the full cover of
sequential integers $1,2,3,..n^{2}$ with RCD sums of $S(n)=n(n^{2}+1)/2$.
From many sources, e.g. the cover of Swetz\cite{Swetz}, and with its
vertical invert $M_{3}$, to which we include the order $3$ addition table, $%
AT_{3}$, of the same elements in which successive rows are augmented by $3$,
all in matrix notation:

\begin{equation}
Luoshu=\left[ 
\begin{array}{ccc}
4 & 9 & 2 \\ 
3 & 5 & 7 \\ 
8 & 1 & 6%
\end{array}%
\right] ,M_{3}=\left[ 
\begin{array}{ccc}
8 & 1 & 6 \\ 
3 & 5 & 7 \\ 
4 & 9 & 2%
\end{array}%
\right] ,AT_{3}=%
\begin{bmatrix}
1 & 2 & 3 \\ 
4 & 5 & 6 \\ 
7 & 8 & 9%
\end{bmatrix}%
,E_{3}=\left[ 
\begin{array}{ccc}
1 & 1 & 1 \\ 
1 & 1 & 1 \\ 
1 & 1 & 1%
\end{array}%
\right] .
\end{equation}%
The first and smallest $Luoshu$ magic square is the sole $3$-by-$3$ magic
square (see Andrews\cite{Andrews}[MSC1,2], Swetz\cite{Swetz5},\cite{Swetz})
and dates before the \textit{Warring States} period in China 403-221 BCE,
and possibly even two millenia earlier. We also added $E_{3}$, a constant
order $3$ matrix of all $1$'s that will soon prove useful.

Both $Luoshu$\ (sometimes called $Lo$ $Shu$) and $M_{3}$ have RCDs of $15$,
but $AT_{3}$ does not for its outer rows and columns, and so is not magic,
but nevertheless affords an example of a pandiagonal\cite{WeissteinPan}
square in which all the continued broken diagonals have the same sum as the
main diagonals. The pandiagonal property is easily seen by placing a copy of 
$AT_{3}$\ to its right:%
\begin{equation}
\begin{array}{ccc}
1 & \mathbf{2} & 3 \\ 
4 & 5 & \mathbf{6} \\ 
7 & 8 & \mathit{9}%
\end{array}%
\begin{array}{ccc}
1 & \mathit{2} & 3 \\ 
\mathit{4} & 5 & 6 \\ 
\mathbf{7} & 8 & 9%
\end{array}%
\end{equation}

and noting that the parallels of the diagonals of $AT_{3}$, i.e. $\mathbf{2}+%
\mathbf{6}+\mathbf{7}$ , $3+\mathit{4}+8$, $1+\mathbf{6}+8$, and $\mathit{2}+%
\mathit{4}+\mathit{9}$ of this rectangle all have the same sum of $15$, the
RCD of $M_{3}$, which is not pandiagonal. In the Appendix we discuss an
order 4 magic square which is pandiagonal and exhibits some of the same
trends under compounding that we are able to demonstrate with our main theme.

$M_{3}$ (and $Luoshu$) are unique aside from their $8$ variants under
rotations and reflections. These are not counted as distinct in most
literature on magic squares and are consistent with most modern literature
in running $1,2,..,3^{2}=9$.

\subsection{Movement in $M_{3}$ (and $Luoshu$)}

Here one notes the Knight's move sequence from $1$ to $2$, followed by
another from $2$ to $3$, then a jump to $4$ before sliding along the
diagonal $4,5,6$ followed by another jump to 7, followed by two more
Knight's moves, $7$ to $8$ and finally $8$ to $9$. The RCDs are all $15$.
Swetz\cite{Swetz} has described the movement from the first cell to the last
as steps in the `Yubu' dance.

\subsection{Associative (or Regular) Magic Squares (AMSs)}

$Luoshu$ and $M_{3}$ are called associative\cite{Weisstein} (or regular)
magic squares (AMS) as their antipodal pairs all sum to a common value, so
that here $10=1+9=2+8=3+7=4+6$.

\section{Frierson's Associative Compound Magic Squares [CMSs] begin at order 
$3^{2}=9$}

The smallest and most famous CMS is taken from Frierson's chapter 5\cite%
{Frierson} in Andrews MSC1\cite{Andrews}, which is also used as Andrews'
Figure 96 on page 44 from his introductory chapter\cite{Andrews}. We label
it $T_{9A}$,\ where\ $T$\ stands for tessellated and\ our terminology
reflects the $A$ggregation of successive groups of nine integers in a tiled
fashion using the pattern of $M_{3}$,\ augmented by successive increments of
nine times $E_{3}$ in the same pattern, so it is clear that it can be
compacted to a $3$-by-$3$ "compact matrix" using multiples of $9E_{3}$ as:

\begin{center}
\begin{equation}
T_{9A}=%
\begin{tabular}{|l|l|l|}
\hline
$M_{3}+63E_{3}$ & $\ \ \ \ \ \ M_{3}$ & $M_{3}+45E_{3}$ \\ \hline
$M_{3}+18E_{3}$ & $M_{3}+36E_{3}$ & $M_{3}+54E_{3}$ \\ \hline
$M_{3}+27E_{3}$ & $M_{3}+72E_{3}$ & $M_{3}+9E_{3}$ \\ \hline
\end{tabular}%
,
\end{equation}
\end{center}

which is one of six related CMS's in three pairs found by Frierson,\textbf{\
which pairs exhibit a different information entropy - }see later discussion
under "Spectral Measures".

It is also an associative magic square, with RCDs of $369$, since all RCDs
sum the linesums of: $3M_{3}+108E_{3}=369$. This associative property is
preserved in all the larger Frierson Compound Magic Squares (CMSs) studied
here that follow for orders of the powers of $3$, i.e. $9,27,81,..$. of
which Frierson found $6$ at order $3\ast 3=9$, and we find $90$ at order $%
3\ast 9=27$, then $2520$ at order $81$, ...

\subsection{$T_{9A}$ in its explicit $9$-by-$9$ form is quite ancient}

\begin{equation}
T_{9A}=%
\begin{tabular}{|l|l|l|}
\hline
\begin{tabular}{lll}
$71$ & $64$ & $69$ \\ 
$66$ & $68$ & $70$ \\ 
$67$ & $72$ & $65$%
\end{tabular}
& $\ \ \ 
\begin{tabular}{lll}
$\mathbf{8}$ & $\mathbf{1}$ & $\mathbf{6}$ \\ 
$\mathbf{3}$ & $\mathbf{5}$ & $\mathbf{7}$ \\ 
$\mathbf{4}$ & $\mathbf{9}$ & $\mathbf{2}$%
\end{tabular}%
$ & 
\begin{tabular}{lll}
$53$ & $46$ & $51$ \\ 
$48$ & $50$ & $52$ \\ 
$49$ & $54$ & $47$%
\end{tabular}
\\ \hline
\begin{tabular}{lll}
$26$ & $19$ & $24$ \\ 
$21$ & $23$ & $25$ \\ 
$22$ & $27$ & $20$%
\end{tabular}
& 
\begin{tabular}{lll}
$44$ & $37$ & $42$ \\ 
$39$ & $41$ & $43$ \\ 
$40$ & $45$ & $38$%
\end{tabular}
& 
\begin{tabular}{lll}
$62$ & $55$ & $60$ \\ 
$57$ & $59$ & $61$ \\ 
$58$ & $63$ & $56$%
\end{tabular}
\\ \hline
$%
\begin{tabular}{lll}
$35$ & $28$ & $33$ \\ 
$30$ & $32$ & $34$ \\ 
$31$ & $36$ & $29$%
\end{tabular}%
$ & 
\begin{tabular}{lll}
$80$ & $73$ & $78$ \\ 
$75$ & $77$ & $79$ \\ 
$76$ & $81$ & $74$%
\end{tabular}
& 
\begin{tabular}{lll}
$17$ & $10$ & $15$ \\ 
$12$ & $14$ & $16$ \\ 
$13$ & $18$ & $11$%
\end{tabular}
\\ \hline
\end{tabular}%
.
\end{equation}

CMSs of multiplicative order $mn$, whose tiled subsquares of orders $m$ and $%
n$ are also magic squares within each subsquare, are found back to at least
the 10th century CE in Persia for the smallest order 9 case $m=n=3$, see
Swetz\cite{Swetz}, and a top-bottom reflection of $T_{9A}$ was recorded in
Arabia by Abul-Wafa al-Buzani (940-997/8 CE), and is found in Descombes\cite%
{Descombes} (p.253/4) and Sesiano \cite{Sesiano}. See also Lam\cite{Lam} and
Li Yen\cite{LiYen}.

A partner CMS to $T_{9A}$ appears in Frierson's 1907 paper, $T_{9D}$, for
which it helps to consider the order $9$ addition table of the first $81$
integers, $AT_{9}$, an obvious generalization of $AT_{3}$ with a first row
of $1,2,..9$, since in 1960 Cammann\cite{Cammann} suggested that this first
CMS may also have been constructed from the rows of such a table.

Judging from the lack of citations to prior work, the first type has been
rediscovered, apparently independently, by many authors, particularly over
the past two centuries, with the partner CMS rarely mentioned. The original
method of construction may have been done by incrementing the upper middle
subsquare by $9$'s and placing them in the same pattern as $M_{3}$. In 2002
Chan and Loly [CL] \cite{Chan} realized that this construction of $T_{9A}$\
is likely why it has been rediscovered many times. The row and column line
sums are clearly magic, as also are the diagonal line sums. CL also
suggested that this amounted to a fractal construction.

Frierson further established a partner square, $T_{9D}$, as well as two
other pairs that we will discuss shortly for a sextet $T_{9A,D,B,E,C,F}$,
while for this order $9$ several estimates indicate an astronomically large
number of $7.8448(38)\times10^{79}$ magic squares at this order for which we
note Walter Trump's website\cite{Trump}.

We note that in 1908 Andrews\cite{Andrews} stated: "The writer believes that
these highly ingenious combinations were first devised by Prof. Hermann
Schubert\cite{Schubert}", whose publication dates to 1898, a sentence which
was deleted in the 1917 edition\cite{Andrews}, but we now know that they
were at least some 900 years earlier! We also found that W.H.Thompson\cite%
{Thompson} constructed the vertical invert of $T_{9A}$ in 1869. Also Pickover%
\cite{Pickover} gave this CMS in 2002, the same year as CL's\cite{Chan}
first CMS report - see Pickover's Chapter 2 (page 81=9*9!) but without
attribution to Frierson, who is only mentioned later on his pages 222,3 for
an order 8 MS! Pickover also used the section title "Composite (or Compound)
Squares", and on the same page says "This reminds me of a fractal,\ ...".

Cammann\cite{Cammann} pointed out that the sums of the subsquares in $T_{9A}$%
, themselves magic squares, also form a magic square, and staring us in the
face are in fact many other CMSs due to rotation and reflection of each $3-$%
by$-3$ subsquare.

\subsection{Frierson's Sextet}

Frierson\cite{Frierson} arranged his six CMSs in vertical pairs on a single
diagram (his page 134):

\begin{center}
\begin{tabular}{|l|l|l|}
\hline
$T_{9A}$ & $T_{9B}$ & $T_{9C}$ \\ \hline
$T_{9D}$ & $T_{9E}$ & $T_{9F}$ \\ \hline
\end{tabular}%
.

Table 1 - Schematic of Frierson's display of the sextet of $F_{9}$'s.
\end{center}

This style will be useful later for order $27$ CMSs. It is easy to see how
the top centre $M_{3}$ subsquare of (4) incremented by $9$'s and a Knight's
move to the lower right subsquare, etc. in the same pattern as $M_{3}$
itself, followed by another increment of $9$ in another Knight's move to the
middle of the LHS. Then down to the bottom LHS with another increment of $9$%
, before incremented diagonal moves up to the to RHS, followed by drop to
the middle RHS, followed by two more Knight's moves to the upper RHS and
finally the bottom centre.

For order $9$, $T_{9A}$ above is the most obvious of CMSs given the
uniqueness of a single order $3$, versus the existence of $880$ distinct
order $4$ magic squares, and thus many more CMSs.

\subsection{Our extensions to larger CMSs}

We discuss here an infinite family of orders $n$ of the powers of three, $%
n=3^{l}$, $l=1,2,..$, as an extension of the 1907/8 study by Frierson, and a
1997 sequel by Bellew \cite{Bellew}, who cited Andrews\cite{Andrews}, 
\textbf{but not Frierson} explicitly!

We are able to give unprecedented insight into a remarkable family of MS's
that began a full millenium ago by extending the spectral properties of
magic squares treated as matrices in LAA\cite{LAA}, and extended by our
later studies \cite{CRL} [CRL] using the more powerful singular value
decomposition (SVD) \cite{CRL}, which demonstrates clearly the SV clan
structure of this family by extending the algebraic formulation begun by
Frierson\cite{Frierson}.

\section{Frierson's parameterization of the order 3 in Andrews \textit{Magic
Squares and Cubes}}

Our starting point is a paper in \textit{The Monist} journal (editor Paul
Carus) in 1907 by Lorraine Screven Frierson\cite{Frierson} of Shreveport\cite%
{Shreveport}, Louisiana, which extends a parameterization of the smallest
magic square of order $3$ in order to generate a sextet of related magic
squares of order $3\ast 3=9$, called Compound (or Composite) Magic Squares
[CMS]. This article was reprinted in the classic \textit{Magic Squares and
Cubes} by W.S.Andrews\cite{Andrews}[MSC1] 1908, with chapters and sections
by W.S.Andrews, L.S.Frierson and others, which are essentially edited
versions of papers originally published in \textit{The Monist} from 1905 for
most of a decade. MSC1 is a critical reference despite shortcomings in
referencing even to \textit{The Monist}, whose included papers share the
same lack of references to earlier sources. A second edition in 1917 will be
denoted MSC2, and was reprinted more recently.

For $l=1,n=3$ we use Frierson's notation \cite{Frierson}\ but omit a common
constant $c$\ in all cells by starting with $0$ at the top centre, followed
by a Knight step to the bottom right corner, placing $y$,\ and a second
Knight step $y$ increment to the middle of the left column for $2y$.
Different increments, $v$ and $2v$ are made on the opposite sides by another
pair of Knight moves from top centre. Then the centre cell is made the
average of those to its left and right, i.e. $v+y$, so that the linesum of
this middle row and column is then $3(v+y)$. The remaining cells are
completed so that all the RCDs have this same linesum. This is shown in
matrix notation to facilitate later spectral function operations:

\begin{equation}
M(v,y)=\left[ 
\begin{array}{ccc}
2v+y & 0 & v+2y \\ 
2y & v+y & 2v \\ 
v & 2v+2y & y%
\end{array}%
\right] ,
\end{equation}

where we will call the pair $v,y$ a "\textbf{couple}". It is clear that
swapping $v,y$ exchanges the first and third columns, i.e. flips the square
from left to right. As already noted, these are not regarded as distinct.

Frierson's parameterization of 3rd order can now be written as:

\begin{equation}
F_{3}(k,v,y)=kE_{3}+M(v,y),  \label{F}
\end{equation}

and since we use $v,y>0$, $k$ gives the smallest entry, which is usually
chosen as either $1$ here (or $0$ by some authors).

[We note that Bellew\cite{Bellew} in 1997 used capital $V,Y$ \ variables
instead of our lowercase $v,y$.]

\textbf{We have used the constant term }$k$\textbf{\ in place of Frierson's }%
$x$\textbf{\ since }$x$\textbf{\ has a standard use in matrix calculus which
is needed later, }and $c$ was used in an alternate parameterization by Lucas%
\cite{Lucas}, described in the Appendix.\textbf{\ }

Then\ $F_{3}(1,3,1)$ is $M_{3}$, while $F_{3}(1,1,3)$\ is the $Luoshu$. Note
that swapping $v$ and $y$ swaps the left and right columns of\textbf{\ }$%
F_{3}$\textbf{. }

Then\textbf{\ }$F_{3}$ has the following properties:

\begin{itemize}
\item it is associative - and this property is maintained in the iterative
compounding of this paper,

\item rotating or flipping $F_{3}$ about the centre still describes all $8$
possible variants of third order magic squares under rotation and reflection,

\item the RCDs all sum to $3k+3v+3y=3(k+v+y)$.
\end{itemize}

The present study extends the powerful parameterized construction of these
by Frierson\cite{Frierson} for order $3^{2}=9$, to higher orders that are
powers of $3$, and first reported by us [LC] at a 2009 conference\cite%
{Loly2009} with an emphasis on counting the number of such magic squares for
orders $n=3^{l},l=1,2,3,...$,\ \textbf{where }$l$\textbf{\ will now be
called "level"}.

\subsection{Frierson's order 9\ parameterization for the smallest (level $%
l=2 $) CMSs}

CMS's of order $n=3^{l}$, $l=2,3..$ are constructed in an iterative manner
from the fundamental parameterized order $3$ pattern, as done in Frierson's 
\cite{Frierson} algebraic study of the smallest CMS's of order\textbf{\ }$%
l=2,n=3^{2}=9$, which consist of six distinct CMS's.

Now we are able to give a complete account of those of the powers of $3$
opened up by Frierson\cite{Frierson} and Bellew\cite{Bellew}, now including
their spectra.

Frierson\ \cite{Frierson} repeated the same associative pattern with two
more parameters, $s,t$, by replacing $v,y$ in $\left( \ref{F}\right) $ for
the couple of level $l=2$:

\begin{equation}
M(s,t)=%
\begin{bmatrix}
2s+t & 0 & s+2t \\ 
2t & s+t & 2s \\ 
s & 2s+2t & t%
\end{bmatrix}
\label{F9}
\end{equation}%
to help in describing $n=9$ associative compound magic squares (CMS9's)
whose elements are then used to provide $s,t$ increments to copies of $%
F_{3}(k,v,y)$ placed in the nine order $m=3$ submatrices tiled to fill a
larger $n=9$ matrix, producing a general ninth order associative compound
magic square, $F_{9}(c,v,y,s,t)$ which has a magic linesum of $9(k+v+y+s+t)$%
. We denote this process by:

\begin{equation}
F_{9}(k,v,y,s,t)=kE_{9}+M(s,t)\tbigotimes M(v,y)
\end{equation}%
in which $\tbigotimes $\ is suggested by the Kronecker product formulation
for CMS's described by Rogers and us \cite{RCL}[RCL].

N.B. Bellew\cite{Bellew} used capital $V,Y,S,T$ in place of Frierson's
lowercase variables. Another useful reference for a broader context of
larger component magic squares is given by Derksen, Eggermont and van den
Essen\cite{Derksen}.

Then a 9-by-9 matrix is constructed from the elements that are the sum of
the components stacked vertically in each cell of Frierson's Figure 228,
which includes a common term, here $k$, for his\ $x$, added to each.

\subsection{Frierson's Order $9$ Sextet}

Frierson generated $6$ distinct numerical $F_{9}$ CMSs by adding a second
couple, $s=27,t=9$, to his first couple, $v=3,y=1$,\ to guarantee full cover
(without gaps or overlap). The explicit algebraic form of $F_{9}$\ is
identical with that of Frierson, aside from a common $x$\ in all his
elements which we have replaced when needed by the constant\textbf{\ }$k$
given earlier.

For $n=9$ at level $l=2$ Frierson simply stated that: `Only six forms may be
made, because, excluding our $k$ [his $x$] whose value is fixed, only six
different couples may be made from the four remaining symbols $v,y,s,t$. '.

These $6$ couples are: three for $v$ with $y,s,t$; then two for $y$ with $%
s,t $; and finally $s$ with $t$. Note that this is still at the algebraic
level before specific parameters are used to produce natural CMSs.

Later we show that the entropy $H$ which decreases from the right column to
the left, is the same for vertical pairs.

\begin{itemize}
\item They are associative by construction, as are the individual tiled
subsquares, i.e. all antipodal pair cells sum to $2(v+y+s+t)$, which is
twice the centre cell.

\item The $T_{9}$'s are $T_{9A}=F_{9}(1,3,1,27,9)$, where the first $1$ is
the constant $k$, then $T_{9D}=F_{9}(1,9,27,1,3)$, then $T_{9B}$,$T_{9E}$
and $T_{9C}$,$T_{9F}$ - see our later Table 3.

\item Since $F_{9}$ is associative by construction, the sextet are also, as
are all the tiled $3$-by-$3$ subsquares.

\item Moreover Table 1 contains all the basic ninth order compound magic
squares, aside from variants due to rotations and reflections of subsquares.

\item The centre cell of $F_{9}$\ is the sum of the $4$ variables, $%
k+v+y+s+t $, and is the average of antipodal pairs, while the bottom centre
cell is always twice that expression, less a $k$.

\item The RCD linesum of $9(k+v+y+s+t)$ summing the values $1,3,9,27$ and
adding $k=1$ is $9\times 40+9=369$.
\end{itemize}

\subsection{Coding using Mathematica\protect\cite{Mathematica}}

First the $vy$ pair, then the $st$ pair, and finally their Kronecker product:

$vy[v\_,y\_]:=\{\{2v+y,0,v+2y\},\{2y,v+y,2v\},\{v,2v+2y,y\}\};$

$st[s\_,t\_]:=\{\{2s+t,0,s+2t\},\{2t,s+t,2s\},\{s,2s+2t,t\}\};$

$f9algebra=KroneckerProduct[st[s,t],vy[v,y]]$

to which one adds the constant matrix $kE_{9}$, where $E_{9}$ is the order $9
$ matrix of all $1$'s.

\subsection{Counting the six order 9 CMS's}

Citing only Andrews' book, but not Frierson's article, Bellew \cite{Bellew}
nevertheless uses Frierson's algebra before giving an argument expressed in
integer values of the parameters that since $k$ is fixed (usually at $0$ or $%
1$), distinct values of the $2$ pairs (couples) of parameters $v,y$ and $s,t$
in which both $v$ and $y$, as well as $s$ and $t$, are interchangeable mean
that there are only: 
\begin{equation}
\left[ (4\times 3)/2\right] \left[ (2\times 1)/2\right] =6  \label{six}
\end{equation}%
unique ways to assign the variables, shown later in Table 3.

\section{Beyond Frierson's $n=9$ sextet to order $27$ (level $l=3$)}

Compounding in a similar fashion to Frierson to order $27$ was suggested
briefly by Bellew\cite{Bellew} in 1997, even though such a large square is
rather unwieldy. In fact we chose $p,q$ variables above after Bellew, but
later noted that he used those for Frierson's $s,t$, so we have followed
Frierson here at order $9$, and then we use $p,q$ for order $27$. Clearly
more parameter pairs can be used for orders $81,243,729,..$ which have\ much
larger CMSs.

Bellew\cite{Bellew} actually considered the counting the magic squares for
two themes, the first reviewing Frierson's parameterization for order $9$
CMSs and briefly suggesting its extension to order $27$, which is developed
fully here, but also a second theme for pandiagonal or Nasiq MSs for orders $%
\geq 5$ which included an order $9$.

With this background, and including spectra not included in most earlier
compounding, we can now proceed to our main theme - to give a complete
account of the generalization of Frierson's scheme to the next order of $%
n=27 $, and later we extend (generalize) this logic for all levels $l$,
before using this powerful formulation to give an algebraic account of the
main spectral function, specifically the singular values for entropic
measures.

\subsection{Order $27$}

The extension to $n=27$ follows similarly with the addition of another pair
of parameters in $M(p,q)$ which has a magic linesum of $3(p+q)$. When this
is compounded with $F_{9}$\ it produces $F_{27}$, which is again
associative, and aside from an overall constant term, describes all possible
compound magic squares with tiled subsquares of orders $3$\ and $9$.

N.B. Our use of\textbf{\ }$p,q$ here for order $27$ CMS's is not the same as
Bellew's\cite{Bellew}\ use of $P,Q$ for his discussion of ninth order.

Now extend (16) to the next compound order of $n=27$:

\begin{equation}
F_{27}(k,v,y,s,t,p,q)=kE_{27}+M(p,q)\tbigotimes [M(s,t)\tbigotimes M(v,y)]
\end{equation}

Since $F_{27}$ is rather large to display explicitly here we continue with
our compact representation for order $27$.

\subsection{The lowest entropy case for order $27$}

For order $27$ the obvious generalization of the lowest entropy order $9$
pair adds a pair $p,q$ with $p=243,q=81$, which has $T_{A}=F_{9A}$ \ in the
top middle order $9$ subsquare with versions incremented by multiples $81$
of an order $9$ with all its elements unity, $E_{9}$, placed in the
corresponding cells in the pattern of $M_{3}$ for a compact representation
of an order $27$ matrix (which otherwise are a challenge to exhibit
explicitly):

\begin{center}
\begin{equation}
F_{27A}=%
\begin{tabular}{|l|l|l|}
\hline
$T_{9A}+7\times81E_{9}$ & $\ \ \ \ \ \ \ T_{9A}$ & $T_{9A}+5\times81E_{9}$
\\ \hline
$T_{9A}+2\times81E_{9}$ & $T_{9A}+4\times81E_{9}$ & $T_{9A}+6\times81E_{9}$
\\ \hline
$T_{9A}+3\times81E_{9}$ & $T_{9A}+8\times81E_{9}$ & $\ \ T_{9A}+81E_{9}$ \\ 
\hline
\end{tabular}%
\end{equation}
\end{center}

where now the multiples of $9E_{3}$ in$\ T_{9A}$\ are now multiples of $%
81E_{9}$. Clearly $F_{27}$'s are both $3$- and $9$-partitioned.

\subsection{Counting the 90 order 27 CMS's at level $l=3$\label%
{firstcouples27}}

From the $6$ parameters $v,y,s,t,p,q$ there are $6!=720$ ways of doing this,
of which some are to be counted as 'basic', while others not. We interpret
the logic of Frierson \cite{Frierson} and Bellew \cite{Bellew} as an
extension of $\left( \ref{six}\right) $ to give $90$ $F_{27}$'s:

\begin{equation}
\left[ (6\times 5)/2\right] \left( \left[ (4\times 3)/2\right] \left[
(2\times 1)/2\right] =6\right) =90\mathbf{.}
\end{equation}%
Here there are six first couples, then four second couples, and finally two
third couples. There are now $15$ distinct `first' couples now multiplied by 
$6$ `second' couples, the number found in $F_{9}$. These are counted as
follows: $5$ for $y=1$, $4$ for $y=3$, $3$ for $y=9$, $2$ for $y=27$, $1$
for $y=81$, for a total of $15$, all multiplied by $6$ from the second
couples.

Having extended Frierson style parameterization for the construction of
order $27$ CMSs, we now turn to spectral measures that give deeper insight
into their properties. To proceed further we need the SVs, $\sigma _{i}$,
for orders $9$ and $27$ which we obtained from Mathematica\cite{Mathematica}
and Maple\cite{Maple} symbolic calculations next.

\section{Matrix Properties - Singular Values (SVs, $\protect\sigma _{i}$)
versus Eigenvalues (EVs, $\protect\lambda _{i}$) for Magic Squares}

Our first foray into the spectra of Frierson's CMSs was presented at a 2009
conference only used EVs, but all our subsequent MS,CMS studies now use the
always positive and declining SVs, whose number of non-zero SVs gives the
rank of the matrix. Few of the sources that we find in the literature
examine matrix spectra, except notably Kirkland and Neumann\cite{Kirkland}
in 1995, drawn to our attention by Adam Rogers c.2005 in connection with
MATLAB\cite{MatLab}, which has a $magic[n]$ function that delivers a single
magic square of odd, even and doubly-even orders. Since then we have
progressively shown how SVs lead to powerful measures for comparing magic
squares, first in 2007 by Loly, Cameron, Schindel and Trump [LCTS]\cite{LAA}%
, then a big leap was made by us in extending Shannon information entropy
measures that Newton and DeSalvo\cite{NDS} used for Sudoku matrices\cite{NDS}
to magic square issues in 2012-13 \cite{CRL} [CRL], and most recently
detailed by Rogers, Cameron and Loly\cite{RCL} [RCL] in 2017. See a standard
text such as Horn and Johnson\cite{HJ} for SV background.

\subsection{Eigenvalues - EVs, $\protect\lambda _{i}$}

First we set the determinant of $M_{3}$ less $x$ times the column vector of
three ones\ equal to zero:

\begin{equation}
Det\left[ 
\begin{array}{ccc}
8 & 1 & 6 \\ 
3 & 5 & 7 \\ 
4 & 9 & 2%
\end{array}%
\right] -x\left[ 
\begin{array}{c}
1 \\ 
1 \\ 
1%
\end{array}%
\right] =0\text{,}
\end{equation}%
for the characteristic polynomial:

\begin{equation}
x^{3}-15x^{2}-24x+360=\allowbreak \left( x-15\right) \left( x^{2}-24\right)
=\allowbreak 0,
\end{equation}

for eigenvalues $\lambda _{i}=15,\pm 2\sqrt{6}$, noting that the effect of
rotation and reflection on $M_{3}$ is to change the imaginary eigenvalues to
real ones in an alternating fashion\cite{LAA}.

Since some larger magic squares have just one non-zero eigenvalue, $\lambda
_{1}$, the RCD linesum, in 2017 Loly, Cameron and Rogers\cite{LCR}[LCR]
concluded that the Singular Values, SVs, $\sigma _{i}$, presented as
positive values declining from the linesum, provided a more useful tool for
assessing magic squares than the eigenvalues, as introduced next.

An invitation to give the lead keynote talk at IWMS2007 allowed Loly and
Cameron to show that matrix eigenvalue analysis of highly singular magic
squares needed to be replaced by the more powerful Singular Value [SV]
analysis. Here the squares of the SVs are the EVs of the product of a matrix
and its transpose, for which we refer to Horn and Johnson\cite{HJ}, to
understand 1EV MSs of order 4 and 5, as reported in 2009 by Loly, Cameron,
Trump and Schindel in LAA\cite{LAA}, and fully by Rogers, Cameron and Loly%
\cite{RCL}[RCL] in 2017.

\subsection{Singular Values - SVs, $\protect\sigma _{i}$}

Now the critical Singular Values (SVs, $\sigma _{i}$) which are always
positive or zero, never complex nor imaginary, and will be the same for both 
$F3\left( 1,3,1\right) $ and $F3\left( 1,1,3\right) $. We began to use the
SVs in LCTS\cite{LAA} 2007/9 when encountering magic squares of orders
greater than three with some vanishing EVs, for which the number of non-zero
SVs gives the matrix rank, $r$.

As an example take the matrix product of $M_{3}$ with its transpose:

\begin{equation}
\left[ 
\begin{array}{ccc}
8 & 1 & 6 \\ 
3 & 5 & 7 \\ 
4 & 9 & 2%
\end{array}%
\right] \left[ 
\begin{array}{ccc}
8 & 3 & 4 \\ 
1 & 5 & 9 \\ 
6 & 7 & 2%
\end{array}%
\right] =\left[ 
\begin{array}{ccc}
101 & 71 & 53 \\ 
71 & 83 & 71 \\ 
53 & 71 & 71%
\end{array}%
\right] ,
\end{equation}

and using $X=\sigma ^{2}$ for the characteristic polynomial, $%
X^{3}-255X^{2}+8556X-29\,340$, which is a cubic equation with the
factorization:

\begin{equation}
(X-15^{2})(X-48)(X-12)=0,
\end{equation}

so that the squares of the SVs, $\sigma _{i}^{2}=15^{2},48=3\times
4^{2},12=3\times 2^{2}$, where their positive square roots are the SVs, $%
\sigma _{i}$ always presented in declining positive values:

\begin{equation}
\sigma _{i}=15,4\sqrt{3},2\sqrt{3},\text{with numerical values }%
15,6.9282..,3.464\,..\text{.}
\end{equation}

See LCTS\cite{LAA} and CRL\cite{CRL} for more on SVs, the latter having
decreasing positive values from the leading SV, $\sigma _{1}=15$, which is
the same as the RCD linesum EV. Also the reverse product is different, but
has the same spectra - a useful feature of SVDs! \textbf{The SVs are also
invariant to rotations and reflections of these (square) matrices.}

With the SVs established, we note that our 2013 study\cite{CRL} showed that
the $880$ of order $4$ have $63$ different singular value "clans" in 2013
(CRL) \cite{CRL}. After Loly gave a talk at McGill later in Summer 2009
noting that we had not found any magic squares with rank less than $3$, Sam
Drury\cite{Drury} proved that MSs have a minimum rank of $3$.

\subsection{The couple $v,y$ for level $l=1$}

In preparation for level $l=2$ for $n=9$ it will be useful to examine this
simplest case as follows. Here the linesum SV: $\sigma _{1}=3(k+v+y)$, and
the pair:

\begin{equation}
\sigma _{2,3}^{2}=3(v\pm y)^{2}
\end{equation}

These are included in our later Table 4.

So for $M_{3}$ and $Luoshu$ when $v,y$ are $3,1$, or vice versa, the
(positive) SVs are $\sigma _{1}=15$, and $\sigma _{2,3}=4\sqrt{3},2\sqrt{3}$%
, as already noted above, for full rank $3$.

By contrast, $AT_{3}$ has singular values: $\left[ 16.\,\allowbreak
848..,1.\,\allowbreak 0684..,0\right] \allowbreak $, and rank $2$ - see
Table 2 later.

The present authors and colleagues have extended earlier studies of singular
values spectra of magic squares to the complete set of the $880$, as well as
to selected higher order magic squares at a 2007 conference [LAA\cite{LAA}
2009], and in greater detail at with a virtual presentation at another in
2012 [DMPS2013].

\subsection{Factorization of the SV characteristic polynomial for $T_{9A}$}

Using $X=\sigma ^{2}$ in the characteristic polynomial which factors as:%
\begin{equation}
(X-369^{2})(X-3\times 108^{2})(X-3\times 54^{2})(X-3\times 12^{2})(X-3\times
6^{2})=0,
\end{equation}

for five non-zero SVs shown later as part of Table 6.

\subsection{Numerical calculations for the SVs}

We have used Mathematica\cite{Mathematica} and Maple\cite{Maple}, including
a subset of the latter in the ScientificWorkplace\cite{SWP} [SWP]TeX system
used for the preparation of this manuscript. The Python\cite{Python}
libraries Numpy and Sympy were also used, and we note that other online
tools for calculating SVs include Keisan\cite{Keisan} and "bluebit"\cite%
{bluebit}.

Now we are able to provide other measures related to Shannon information
entropy, which measures the degree of order in a system, and can now show an
asymptotic behaviour for increasingly large order CMSs in the present study
of Frierson's partner CMSs.

\subsection{Spectral Measures - Entropy $H$ and Compression $C$}

With the SVs, $\sigma _{i}$, we can proceed further we introduce some
measures introduced for us in 2010 by Newton and DeSalvo\cite{NDS} [NDS] who
considered Sudoku matrices, which are special order $9$ Latin squares of
elements $1..9$ in every row and column arranged so that each occurs in
every tiled $3$-by-$3$ subsquare. In 2013 we extended NDS to magic squares
of orders $3,4,..9$ as well as Latin squares from orders $2,3,4,5,8,9$ in CRL%
\cite{CRL}.

These powerful measures for assessing different magical squares, notably the
Shannon information entropy $H$\ and a very useful percentage Compression $C$%
,\ which we found in 2010 in NDS for completed Sudoku puzzles which may be
regarded as compounded order three Latin squares (Sudoku appeared in
newspapers c.2004). Useful measures of these matrices are now shown in a
tabular report, whose components will now be defined.

First the SVs, $\sigma _{i}$, are normalized by their sum:%
\begin{equation}
\hat{\sigma}_{i}=\dfrac{\sigma _{i}}{\Sigma _{i}^{n}\sigma _{i}},
\end{equation}

then the Shannon information entropy, $H$, is calculated:

\begin{equation}
H=-\Sigma _{i}^{n}\hat{\sigma}_{i}\ln (\hat{\sigma}_{i}),
\end{equation}

named after Boltzmann's $H$-theorem, and finally a very useful percentage
compression measure:

\begin{equation}
C=(1-\dfrac{H}{\ln (n)})\times 100\%,
\end{equation}

which being bounded between $0\%$\ and $100\%$ is very useful for
comparisons between different magic squares.

See the Appendix for a sample numerical calculation for $M_{3}$.

\subsection{Additional measures $R,L$}

In CRL\cite{CRL} we introduced some integer measures for integer square
matrices based on the sums of the even powers of the SVs, ,

\begin{equation}
L=\Sigma _{i}^{n}\sigma _{i}^{4},
\end{equation}%
and especially its shorter version: 
\begin{equation}
R=\Sigma _{i=2}^{n}\sigma _{i}^{4}=L-\sigma _{1}^{4},
\end{equation}

which is also integer for MSs. These are included in Table 2 below. CRL\cite%
{CRL} called the distinct sets of SVs "clans"\cite{CRL}, which usually have
a distinct value of $R$, except so far only for one pair at order 4.

\subsection{Matrix rank of CMSs}

Drury\cite{Drury} showed that magic squares have a minimum rank of $3$, and
therefore if less than their order $n$, are singular with one or more zero
eigenvalues. In 2017 a theorem was given by Adam Rogers and the present
authors \cite{RCL} [RCL], for understanding the matrix rank of CMSs of
combinations of all orders which gives their rank as the sum of their
component ranks, here for $n=9$\ each $3$ less $1$\ for rank $r=3+3-1=5$.

Our 2009 conference report on Frierson's compound squares\cite{Loly2009}
occured before we encountered the Shannon entropy measures later in 2010, so
did not include these powerful measures for the entropy and compression,
which we later encountered later from 2010 paper by Newton and DeSalvo\cite%
{NDS}. These were \ then used in a conference in 2012 with Adam Rogers in
2013 [RCL]\cite{CRL}, which included a table for order $9$ magic squares
including both $T_{9A}$\ and $T_{9D}$, but without further elaboration.

\section{Matrix Properties for $n=3$}

Our first tabular presentation of the matrix spectra:

\begin{center}
\begin{tabular}{|l|l|l|l|}
\hline
matrix & $5E_{3}$ & $M_{3},Luoshu$ & $AT_{3}$ \\ \hline
$\lambda _{i}$ & $15,0,0$ & $15,\pm 2i\sqrt{6}$ & $\frac{3}{2}(5\pm \sqrt{33}%
)$ \\ \hline
$\sigma _{i}^{2}$ & $225,0,0$ & $225,48,12$ & $\frac{3}{2}(95\pm \sqrt{8881}%
) $ \\ \hline
$\sigma _{1}$ & $15$ & $15$ & $16.8481..$ \\ \hline
$\sigma _{2}$ & $0$ & $4\sqrt{3}=6.928..$ & $1.06837..$ \\ \hline
$\sigma _{3}$ & $0$ & $2\sqrt{3}=3.464..$ & $0$ \\ \hline
$H$ & $0.0$ & $0.937098..$ & $0.22595..$ \\ \hline
$C$ & $100\%$ & $14.7017..\%$ & $79.4332..\%$ \\ \hline
rank, $r$ & $1$ & $3$ & $2$ \\ \hline
$R$ & $0$ & $2448$ & $\allowbreak 1.\,\allowbreak 302\,82..$ \\ \hline
$L$ & $50,625$ & $53,073$ & $80,577$ \\ \hline
\end{tabular}

Table 2 - Matrix properties for $5E_{3}$, $M_{3}$ and $AT_{3}$.

N.B. For $AT_{3}$,\ since $\sigma _{1}$ is not integer then nor is $R$.
\end{center}

For $M_{3}$\ the pair $\sigma _{2}$,$\sigma _{3}$ differ by a factor of $2$,
a feature found in later pairs in Table 6. The $14.7017..\%$ compression for 
$M_{3},Luoshu$ is one of the smallest that CRL\cite{CRL} found in a wide
ranging study of magic squares and Latin squares, while we will see that the
larger CMSs here trend to much higher $C\%$'s than we found for the smallest
CMSs of order $9$ of $48.57..\%$ that we showed earlier\cite{CRL}. Extended
in the present work to orders $27,81,243,...$, we find systematically larger
values that tend towards the uniformity of $100\%$. Since any uniform square
matrix of all $1$'s has full compression of $100\%$, a low compression
reflects a more "lumpy" matrix! \ Most other (larger) magic squares have a
much higher compression, especially the compound magic squares studied here.

N.B. After this table we drop further discusson of the EVs ($\lambda _{i}$)
since the SVs ($\sigma _{i}$) give us all the information needed for the
entropy and compression.

Also all versions of\ $E_{n}$\ have just the linesum EV ($\lambda _{1}$) and
SV ($\sigma _{1}$), both $n$.

It is worth noting that in DMPS we did find lower compressions than for $%
M_{3}$'s $14.7017..\%$ for some order $5$ and $9$ MSs, and that $AT_{3}$'s\
high compression shows the high ordering of its elements, only surpassed by
the completely ordered matrices of identical elements, e.g. the $100\%$ of $%
E_{3}$.

\subsection{Zero-based MSs}

If the elements of a MS are chosen to run $0,1,2,3,..(n^{2}-1)$\ instead of
the $1,2,3,..n^{2}$ used here, then the $\hat{\sigma}_{i}$\ will be smaller
since the RCD's are smaller, so that the entropy will be larger and the
compression smaller, e.g. for $M_{3}^{\prime }$ these change to: $%
H^{^{\prime }}=0.985975,C^{^{\prime }}=10.2527\%$.

\section{Frierson's partner CMSs}

Now fill a new $D_{3}$\ in place of $M_{3}$\ with the elements of $AT_{9}$'s
first column, $1,10,19,28,37,46,55,64,72$, in the $M_{3}$ pattern for $D_{3}$%
, enhanced by simple multiples of $E_{3}$ to obtain $T_{9D}$, a spectral
partner MS to $T_{9A}$:

\begin{equation}
D_{3}=\left[ 
\begin{array}{ccc}
64 & 1 & 46 \\ 
19 & 37 & 55 \\ 
28 & 73 & 10%
\end{array}%
\right] ,T_{9D}=%
\begin{tabular}{|l|l|l|}
\hline
$D_{3}+7E_{3}$ & $\ \ \ \ \ \ D_{3}$ & $D_{3}+5E_{3}$ \\ \hline
$D_{3}+2E_{3}$ & $D_{3}+4E_{3}$ & $D_{3}+6E_{3}$ \\ \hline
$D_{3}+3E_{3}$ & $D_{3}+8E_{3}$ & $D_{3}+E_{3}$ \\ \hline
\end{tabular}%
\text{,}
\end{equation}

and is magic, having the same SVs as $T_{9A}$ - see Table 3 below.

Here we used $D$ to indicate that the elements of the subsquares of $T_{A}$\
have been $D$ispersed to other subsquares in a systematic way.

$T_{9D}$ also has an early date before 1000 CE - Cammann\cite{Cammann} noted
that this magic square was found in China by the 13th CE by Yang Hui 1275
CE, and suggested that $T_{9A}$ and $T_{9D}$ were originally derived from
the order $9$ addition table, $AT_{9}$. See also Table 3 below for its
spectra.

\subsection{Frierson's second pair $T_{9B}=F_{9}(1,27,9,3)$ and $T_{9E}$}

$T_{9B}$ uses the first rows of the left hand subsquares of $AT_{9}$, $1,2,3$
with $28,29,30$ and $55,56,57$, to fill a $B_{3}$ with the $M_{3}$ pattern:

\begin{equation}
B_{3}=\left[ 
\begin{array}{ccc}
56 & 1 & 30 \\ 
3 & 29 & 55 \\ 
28 & 57 & 2%
\end{array}%
\right] ,T_{9B}=%
\begin{tabular}{|l|l|l|}
\hline
$B_{3}+21E_{3}$ & $\ \ \ \ \ \ B_{3}$ & $B_{3}+15E_{3}$ \\ \hline
$B_{3}+6E_{3}$ & $B_{3}+12E_{3}$ & $B_{3}+18E_{3}$ \\ \hline
$B_{3}+9E_{3}$ & $B_{3}+24E_{3}$ & $B_{3}+3E_{3}$ \\ \hline
\end{tabular}%
.
\end{equation}

with SVs: $369,145.\,\allowbreak 49..,135.\,\allowbreak
10..,62.\,\allowbreak 354..,31.\,\allowbreak 177$.., and four zeros, as does
its partner $T_{9E}$, not shown.

\subsection{Frierson's third pair $T_{9C}=F_{9}(1,9,27,3)$ and $T_{9F}$}

$T_{9C}$ then uses $1,2,3$ with $10,11,12$ and $19,20,21$ from the top left
order $3$ subsquare of $AT_{9}$ arranged in the$\ M_{3}$\ pattern:

\begin{equation}
C_{3}=\left[ 
\begin{array}{ccc}
20 & 1 & 12 \\ 
3 & 11 & 19 \\ 
10 & 21 & 2%
\end{array}%
\right] ,T_{9C}=%
\begin{tabular}{|l|l|l|}
\hline
$C_{3}+57E_{3}$ & $\ \ \ \ \ \ C_{3}$ & $C_{3}+33E_{3}$ \\ \hline
$C_{3}+6E_{3}$ & $C_{3}+30E_{3}$ & $C_{3}+54E_{3}$ \\ \hline
$C_{3}+27E_{3}$ & $C_{3}+60E_{3}$ & $C_{3}+3E_{3}$ \\ \hline
\end{tabular}%
.
\end{equation}

Now $T_{9C}=F_{9}(1,9,1,27,3)$, with a "partner" $T_{9F}$ .. see also Table
3, with singular values: $369.0,155.\,\allowbreak 88,124.\,\allowbreak
71,51.\,\allowbreak 962,41.\,\allowbreak 5969$, and four zeros, as does its
partner $T_{9F}$, not shown.

\subsection{Matrix Properties for Frierson's 6 natural $9$th order 'basic' $%
F_{9}$'s $l=2,n=9$}

The next table gives the properties for Frierson's six squares (ordered $v>y$
from $M(v,y)$ where $v,y$ in the second row, and $s>t$ from $M(s,t)$ where $%
s $ is written to the left of $t$), showing pairs of isentropic variants:

\begin{center}
\begin{tabular}{|l|l|l|l|l|l|l|}
\hline
$F_{9}$ & $T_{9A}$ & $T_{9D}$ & $T_{9B}$ & $T_{9E}$ & $T_{9C}$ & $T_{9F}$ \\ 
\hline
$v,y$ & $3,1$ & $27,9$ & $27,1$ & $9,3$ & $9,1$ & $27,3$ \\ \hline
$s,t$ & $27,9$ & $3,1$ & $9,3$ & $27,1$ & $27,3$ & $9,1$ \\ \hline
$C$ & $48.572..\%$ & \TEXTsymbol{<}- & $40.0241..\%$ & \TEXTsymbol{<}- & $%
39.8296..\%$ & \TEXTsymbol{<}- \\ \hline
$H$ & $1.12999..$ & \TEXTsymbol{<}- & $1.31781..$ & \TEXTsymbol{<}- & $%
1.32208..$ & \TEXTsymbol{<}- \\ \hline
$R$ & ${\small 1,301,165,856}$ & {\small \TEXTsymbol{<}-} & ${\small %
797,281,056}$ & {\small \TEXTsymbol{<}-} & ${\small 842,630,688}$ & 
\TEXTsymbol{<}- \\ \hline
\end{tabular}

Table 3 - Matrix Properties for Frierson's $F_{9}$ sextet with RCDs, $%
\lambda _{1},\sigma _{1}=369$.
\end{center}

These $R$ values for $T_{9A,D}$\ agree with our 2017 RCL \cite{RCL}, but
since $R$ becomes much larger for $n=27$,$81$,... it will be dropped
henceforth, with an emphasis on the \% Compression which is always bounded
between $0\%$\ and $100\%$.

Here there are 3 sets of SVs, each with different entropies and
compressions. We interpret the reduced compression and higher entropy values
to show that the order decreases from $T_{9A,D}$, through $T_{9B,E,}$, to $%
T_{9C,F}$ are not quite as ordered as $T_{9A,D}$, but are closer to each
other. Clearly the spectral properties are not changed by swapping the
parameters values of the pairs $y,v$\ and $s,y.$

This gives a deeper insight into Frierson's construction than possible
without the spectra.

Our spectral measures for $T_{9A,B,C,D,E,F}$\ differ from $M_{3}$, with $C\%$
of $14.7\%$,\ having a much greater Compression, almost halfway to the $%
100\% $ of a uniform matrix, e.g. $E_{9}$, a trend that increases as we
explore order $27,81,243,..$ compounding later by continuation of the
fractal pattern underlying this particular system, and apparently becomes
asymptotic at about $1.1677038..$ in our later Table 6.

Next we extend Frierson's ideas to the next level, $l=3$ for $n=27$.

\section{Comparing Spectral Algebras for\ $l=1,2,3$ (or $n=3,9,27$)}

\textbf{We followed Frierson in the use of }$v,y$\textbf{\ and then his }$%
s,t $\textbf{, for order }$9$, \textbf{whereas Bellew\cite{Bellew} used }$%
p,q $\textbf{\ instead of Frierson's }$s,t$\textbf{, so we now use }$p,q$%
\textbf{\ for the step to order }$27$\textbf{.}

It is clear that this process could be continued for orders $81$, $243$, ...
but already a clear pattern has emerged which renders that unnecessary as
the next Table will show!

On the basis of Maple and Mathematica calculations we can now state the
formulae for the singular values $\left( n=3^{l}\right) $ of all orders of
Frierson compound squares which consists of the linesum eigenvalue, and $l$
signed pairs and rank:

\begin{center}
\begin{tabular}{|l|l|l|l|}
\hline
$l$ & $1$ & $2$ & $3$ \\ \hline
$n=3^{l}$ & $3$ & $9$ & $27$ \\ \hline
$r=2l+1$ & $3$ & $5$ & $7$ \\ \hline
$S(n)$ & $15$ & $369$ & $9855$ \\ \hline
$\sigma _{1}-nk$ & $3(v+y)$ & $9(v+y+s+t)$ & $27(v+y+s+t+p+q)$ \\ \hline
$\sigma _{2,3}^{2}$ & $3(v\pm y)^{2}$ & $27(v\pm y)^{2}$ & $243(v\pm y)^{2}$
\\ \hline
$\sigma _{4,5}^{2}$ &  & $27(s\pm t)^{2}$ & $243(s\pm t)^{2}$ \\ \hline
$\sigma _{6,7}^{2}$ &  &  & $243(p\pm q)^{2}$ \\ \hline
\end{tabular}

Table 4 -Singular Values for $n=3,9,27$.
\end{center}

In this table the $\sigma _{i}$ for $i>1$ increase by a factor of $3$, so
that their squares increase by factors of $9$. It is clear how this table
can be extended by adding extra pairs, e.g. $a,b;c,d$, etc. for $l=4,5,..$
Considering orders $3,9,27$ in Table 6 above where it does not matter for
the SVs if $v$ is greater or less than $y$ (because of the squares in the
formulae for $\sigma _{2,3}^{2}=3(v\pm y)^{2}$), nor similarly their
positive numerical magnitudes.

N.B. While numerical data for the SVs are usually listed in descending
magnitude the magnitudes of $p,q,s,t,v,y$\ vary in the next table the
magnitudes of $\lambda _{6,7},\sigma _{6,7}$, $\lambda _{_{4,5\text{,}%
}},\sigma _{4,5\text{,}}$ $\lambda _{2,3},\sigma _{2,3}$ will rarely be
sequential!

\subsection{Numerical $F_{27}$ spectra}

Calculations were done with SV formulae in Table 4 above. All have rank $%
7=(3+3-1=5)+3-1$ in agreement with RCL\cite{RCL}.

Case A has the lowest entropy and its counterparts for different orders will
be our main focus. Note that this case has two sets of parameters, $v,s,p$
and $y,t,q$, increasing monotonically.

The integer index $R$ devised by Loly in CRL\cite{CRL} as the sum of the 4th
powers of the SVs (less the one for the linesum) is rather long and we note
just the two extremes:

$R($A$)$= $691,492,899,739,824$ with ln[$R($A$)$]=$34.169874...$,

and $R($O$)$=$420,327,995,019,696$ with ln[$R($O$)$]=$33.672056...$,

from which we conclude that $H$ and especially $C\%$ are more useful in
comparing large MSs than the huge integer $R$'s!

The $90$ order $27$'s would need a $6$-by-$15$ table of $6$ rows for the
isentropic squares and $15$ columns of the different entropies which are now
listed:

\begin{center}
\begin{tabular}{|l|}
\hline
$%
\begin{tabular}{lllllllll}
& $v$ & $y$ & $s$ & $t$ & $p$ & $q$ & $H$ & $C\%$ \\ 
A & $1$ & $3$ & $9$ & $27$ & $81$ & $243$ & $1.16247$ & $64.7291$ \\ 
B & $1$ & $27$ & $3$ & $9$ & $81$ & $243$ & $1.20646$ & $63.3944$ \\ 
C & $1$ & $9$ & $3$ & $27$ & $81$ & $243$ & $1.20697$ & $63.3788$ \\ 
D & $1$ & $3$ & $81$ & $27$ & $243$ & $9$ & $1.34763$ & $59.1110$ \\ 
E & $1$ & $3$ & $243$ & $27$ & $81$ & $9$ & $1.35191$ & $58.9813$ \\ 
F & $1$ & $243$ & $9$ & $3$ & $81$ & $27$ & $1.38498$ & $57.9778$ \\ 
G & $1$ & $9$ & $243$ & $3$ & $81$ & $27$ & $1.38566$ & $57.9573$ \\ 
H & $1$ & $81$ & $3$ & $9$ & $27$ & $243$ & $1.38973$ & $57.8338$ \\ 
I & $1$ & $9$ & $81$ & $3$ & $243$ & $27$ & $1.39035$ & $57.8149$ \\ 
J & $1$ & $243$ & $81$ & $3$ & $9$ & $27$ & $1.46991$ & $55.4010$ \\ 
K & $1$ & $81$ & $27$ & $9$ & $243$ & $3$ & $1.46996$ & $55.3995$ \\ 
L & $1$ & $243$ & $3$ & $27$ & $9$ & $81$ & $1.47129$ & $55.3593$ \\ 
M & $1$ & $27$ & $81$ & $9$ & $243$ & $3$ & $1.47148$ & $55.3533$ \\ 
N & $1$ & $81$ & $3$ & $27$ & $9$ & $243$ & $1.47178$ & $55.3443$ \\ 
O & $1$ & $27$ & $81$ & $3$ & $243$ & $9$ & $1.47193$ & $55.3398$%
\end{tabular}%
$ \\ \hline
\end{tabular}

Table 5 - $F_{27}$ Matrix spectral measures for $15$ clans\ with the lowest
entropy at top and highest at bottom.
\end{center}

In an Appendix Browne's$\cite{Browne}$ order $27$ is shown to have $%
v=27,y=1;s=3,t=81;p=9,q=243$, so it is a variant of case "O" with the
highest entropy, one of $90/6=15$ variants - see the next section.

\subsection{Collecting the lowest entropy sets for higher values of $l,n$}

Late in 2019 we realized that the SVs of higher order versions of the 
\textbf{lowest entropy members}, e.g. $T_{9A}$, $T_{27A}$, .., could be used
directly to obtain the Compression and entropy values so now the SVs for
each pair differ by the same factor of $2$ found in Table 2 for $l=1$, and
these SVs increase by a factor of $27$ as $l$ increases, while the SV's of
each higher pair increase by a factor of $9$ for every increase in $l$. Here
we see these lowest entropies slowly increasing with order $n$\ from $%
0.937.. $\ to $1.168..$\ and clearly becoming asymptotic - a feature that we
now see was probably present in our earlier CRL study\cite{CRL} for magic
squares obtained form the MATLAB's $magic[n]$ function\cite{MatLab}, where
its "Figure 1" showed a slowing increase of the entropies of odd order to $%
n=99$\ from $H=0.937..$\ to $\symbol{126}3.5$\ (which only included the sole
order $M_{3}$ in the present study since those for $n=9,27,..$\ lie well
above our lowest entropy members: for $n=9$ c. $1.8$.., for $n=27$ c. $2.7..$
and for $n=81$.c. $3.22..${\small ).}

\begin{center}
\begin{tabular}{|l|l|l|l|l|l|}
\hline
${\small l}$ & ${\small 1}$ & ${\small 2}$ & ${\small 3}$ & ${\small 4}$ & $%
{\small 5}$ \\ \hline
& $M_{3}$ & $F_{9A}$ & $F_{27A}$ & $F_{81A}$ & $F_{243A}$ \\ \hline
${\small n=3}^{l}$ & ${\small 3}$ & ${\small 9}$ & ${\small 27}$ & ${\small %
81}$ & ${\small 243}$ \\ \hline
{\small RCD} ${\small \sigma }_{1}$ & ${\small 15}$ & ${\small 369}$ & $%
{\small 9855}$ & ${\small 265,761}$ & ${\small 7174575}$ \\ \hline
${\small \sigma }_{2}/\sqrt{3}$ & ${\small 4}$ & ${\small 108}$ & ${\small %
2916}$ & ${\small 78732}$ & ${\small 2125764}$ \\ \hline
${\small \sigma }_{3}/\sqrt{3}$ & ${\small 2}$ & ${\small 54}$ & ${\small %
1458}$ & ${\small 39366}$ & ${\small 1062882}$ \\ \hline
${\small \sigma }_{4}/\sqrt{3}$ &  & ${\small 12}$ & ${\small 324}$ & $%
{\small 8748}$ & ${\small 236196}$ \\ \hline
${\small \sigma }_{5}/\sqrt{3}$ &  & ${\small 6}$ & ${\small 162}$ & $%
{\small 4374}$ & ${\small 118098}$ \\ \hline
${\small \sigma }_{6}/\sqrt{3}$ &  &  & ${\small 36}$ & ${\small 972}$ & $%
{\small 26244}$ \\ \hline
${\small \sigma }_{7}/\sqrt{3}$ &  &  & ${\small 18}$ & ${\small 486}$ & $%
{\small 13122}$ \\ \hline
${\small \sigma }_{8}/\sqrt{3}$ &  &  &  & ${\small 108}$ & ${\small 2916}$
\\ \hline
${\small \sigma }_{9}/\sqrt{3}$ &  &  &  & ${\small 54}$ & ${\small 1458}$
\\ \hline
${\small \sigma }_{10}/\sqrt{3}$ &  &  &  &  & ${\small 324}$ \\ \hline
${\small \sigma }_{11}/\sqrt{3}$ &  &  &  &  & ${\small 162}$ \\ \hline
${\small \sigma }_{total}$ & $26.3923..$ & $680.76..$ & $18366.3..$ & $%
4.9847..10^{4}$ & $1.3387..10^{7}$ \\ \hline
${\small C\%}$ & ${\small 14.7017..}$ & ${\small 48.572..}$ & ${\small %
64.7291..}$ & ${\small 73.4364..}$ & ${\small 78.7368..}$ \\ \hline
${\small H}$ & ${\small 0.93709..}$ & ${\small 1.1299..}$ & ${\small %
1.16247..}$ & ${\small 1.16732..}$ & ${\small 1.1677038..}$ \\ \hline
${\small r=2l+1}$ & ${\small 3}$ & ${\small 5}$ & ${\small 7}$ & ${\small 9}$
& ${\small 11}$ \\ \hline
\end{tabular}

Table 6 - The lowest entropy members of Frierson-type CMSs.
\end{center}

Our main goal in going beyond Frierson's order $9$ CMSs to a full account of
order $27$ is now complete, but we can now make a further extension for the
lowest entropy (highest order) cases.

Since this now completes $n=27$, we will now extrapolate to higher orders -
see later for Sloane's\textbf{\ }A000680\cite{Sloane680} and counting the
isentropic variants illustrated here for $n=9,27$..

\subsection{Asymptotic behaviour}

For $F_{729A}$ with $l=6,n=729$ we find ${\small C\%=82.2829..}$, ${\small %
H=1.167856..}$. The entropy $H$ is clearly flattening out to about $1.168..$%
, while the Compression $C\%$ continues to increase more slowly towards $%
100\%$.

Other CMSs using $T_{B,C,E,F}$, which begin with higher values of entropy,\
compounded with or without $T_{A,D}$, are expected to generate larger
entropies than found above and are not pursued here.

\section{Counting for $n=3^{l}$}

The number of $F_{n}$'s at level $l$\ is the product of the number of first
couples at level $l$, column 3 in the table below, the number of $F_{n}$%
{\small 's} at the previous level ($l-1$), for $l\left( 2l-1\right) $ first
couples, as shown in column $4$, and the number of distinct SV sets in
column 5:

\begin{center}
\begin{tabular}{|l|l|l|l|l|}
\hline
${\small n}$ & ${\small l}$ & $1${\small st couples} & {\small number of }$%
F_{n}${\small 's} & {\small no. of SV sets} \\ \hline
&  & $l\left( 2l-1\right) $ & $\left( 2l\right) !/2^{l}$ & $\left(
2l-1\right) !!$ \\ \hline
${\small 3}$ & ${\small 1}$ & ${\small 1}$ & ${\small 1}$ & ${\small 1}$ \\ 
\hline
${\small 9}$ & ${\small 2}$ & ${\small 6}$ & ${\small 6}$ & ${\small 3}$ \\ 
\hline
${\small 27}$ & ${\small 3}$ & ${\small 15}$ & ${\small 15\times 6=90}$ & $%
{\small 5\times 3=15}$ \\ \hline
${\small 81}$ & ${\small 4}$ & ${\small 28}$ & $28\times 90={\small 2520}$ & 
$7\times 15={\small 105}$ \\ \hline
${\small 243}$ & ${\small 5}$ & ${\small 45}$ & $45\times {\small 2520=113400%
}$ & $9\times 105={\small 945}$ \\ \hline
\end{tabular}

Table 7 - Counting couples, $F_{n}$'s and SV sets{\small .}

Note that $15$ in columns $3$ and $5$ is a coincidence.
\end{center}

Also only the $n=3,9,27$ results in columns 4,5 have been verified, and
those prompted the formulae above and "OEIS" described next.

\subsection{Integer Sequences - we use order $n=3^{l}$ in this paper}

`The On-Line Encyclopedia of Integer Sequences', "OEIS", gives the following
information on the three integer sequences used here:

\subsubsection{Counting $1$st couples}

Sloane's\textbf{\ }\cite{Sloane384}\textbf{\ }A000384: \textit{Hexagonal
numbers}: $n(2n-1)$:

$0,1,6,15,28,45,..$,

in our notation: $l(2l-1)$, and ignoring the zero!

\subsubsection{Counting the number of $F_{n}$'s}

Sloane's \cite{Sloane680}\textbf{\ }A000680: $(2n)!/2^{n}$:

$1,1,6,90,2520,113400,7484400,681080400,81729648000,..$,

in our notation: $\left( 2l\right) !/2^{l}$, again ignoring the first '$1$'.

\subsubsection{Number of SV sets}

Sloane's\cite{Sloane1147} A001147: \textit{Double factorial of odd numbers}:

$a(n)=(2n-1)!!=1\times3\times5\times...\times(2\ast n-1)$:

$1,1,3,15,105,945,10395,135135,..$,

and again ignoring the first '$1$'.

\subsection{Factors of $8$ for $F_{9}$'s, $F_{27}$'s ...}

Bellew's factors of $8$ drew our attention\cite{Loly2009} to the
significance of this aspect of compounding.

Now we note the effect of rotations and reflections of each subsquare, $m=3$
for $F_{9}$'s for $8^{9}$ variations in $F_{9}$'s,\ and both $m=3,9$
subsquares for $F_{27}$'s which now give a factor of $8^{81+9}=8^{90}$
variants of each basic $F_{27}$ due to a factor of $8$ for each of the $9$ $%
m=9$\ subsquares multiplying the factor from $81$ $\ m=3$ subsquares.

Here for $n=27$ we have resolved disparate counts of $8^{18}$ of Trigg \cite%
{Trigg} (1980) and Bellew of $8^{81}$\ to a new result of $8^{81+9}=8^{90}$
by taking account of all orders of tiled subsquares, before generalizing
this for all $l$.

Then for $F_{81}$'s we predict an additional factor of $8^{729}$ for a total 
$8^{729+81}=8^{819}$. We observe that the exponents $9,81,819,..$ may be
found in Sloane's\cite{Sloane523386}

\subsubsection{Number of variants due to subsquare rotations and reflections}

A0523386: \textit{Number of integers from} $1$ to $10^{n}-1$ that lack $0$
as a digit: $0,9,90,819,7380,66429,597870,..$ (ignoring the initial zero).

We also expect that these rotations and reflections of the magic subsquares
in $F_{9}$ will increase the rank of the resultant CMS variants.

\section{ CMSs and\ Fractal patterns c.2000}

Earlier Chan and Loly\cite{Chan} [CL] revived the compounding idea by using
a pandiagonal order $4$ and Euler's 1779 pandiagonal order $7$ to produce an
aggregated CMS of order $12,544=4^{4}\ast 7^{2}$, suggesting that this
process is fractal\cite{Mandelbrot}, i.e. self-similar on all scales, in
order to break records for large magic squares. CL also gave an argument for
the preservation of pandiagonality on compounding that parallels our present
observation of the preservation of associativity on compounding, and while
referencing the important 1997 work of Bellew\cite{Bellew}, focussed on his
treatment of pandiagonal magic squares (PMSs), defined later, rather than
Frierson's associative squares of concern here. They referenced a then
recent paper 1997 paper by Bellew\cite{Bellew} as well as Andrews\cite%
{Andrews}, neither of which explicitly referenced Frierson's parametric
compounding of the order $3$ to order $9$.

\textbf{THIS COMPLETES\ OUR\ EXTENSIONS OF FRIERSON'S and BELLEW's IDEAS.}

\section{CONCLUSION}

Frierson's parameterization set the stage for our generalization here.
Extending his algebraic formulation from order $9$ to highet powers of $3$
has enabled us to project asymptotic behaviour for the lowest entropy
members of this infinite family of CMSs of orders\textbf{\ }$3^{l}$\textbf{, 
}giving the first full account of order $27$..

Our present achievement may be considered somewhat parallel to Ollerenshaw
and Br\'{e}e's\cite{OllyBree} comprehensive study of Most-Perfect
Pandiagonal [MPPD] MSs of orders all multiples of $4$, but enhanced here
with an account of the spectral properties\textbf{. }A preliminary study of
compounding of one of those at order $4$ in our Appendix indicates similar
asymptotic behaviour, which suggests a new look at their parameterization
would be valuable, so the present work will be followed by a study of
parameterizing order $4$ MSs by Ian Cameron\cite{Cameron} using the 1910
scheme of Bergholt\cite{Bergholt}.

\section{Acknowledgements}

We thank Adam Rogers\cite{RCL} [RCL] and Wayne Chan\cite{Chan} [CL] for
earlier compound collaborations. PDL has seen a copy of the original notes
of the present topic and other L.S.Frierson material at Shreveport,
Louisiana, courtesy of Fermand M. Garlington II, \textit{Archives and
Special Collections}, Louisiana State University in Shreveport. "Frierson"%
\cite{Shreveport} is also an unincorporated community and Census-Designated
Place (CDP) in DeSoto Parish, Louisiana, United States.

Loly also received early encouragement from John Hendricks\cite{Hendricks},
originally from our city of Winnipeg, and Harvey Heinz\cite{Heinz}, who were
early active members of a large online recreational mathematics community, \
whom we hope to encourage to include our spectral measures in their future
investigations, including several extensive websites, e.g Harry White\cite%
{White} and Walter Trump\cite{Trump}.

Email: loly@umanitoba.ca, Ian.Cameron@umr.umanitoba.ca.

\appendix

\section{An earlier parameterization by \'{E}douard Lucas in 1894}

Another parameterization for order $3$ by Lucas\cite{Lucas} should be noted
and was drawn to our attention in detail by Sallows\cite{Sallows}, who used
a parameter $c$ which plays the role of our $k$ and Frierson's $x$:

\begin{center}
\begin{equation}
Lucas(a,b,c)=%
\begin{tabular}{|l|l|l|}
\hline
$\ \ c+a$ & $c-a-b$ & $\ \ c+b$ \\ \hline
$c-a+b$ & $\ \ \ \ \ \ c$ & $c+a-b$ \\ \hline
$\ \ c-b$ & $c+a+b$ & $\ \ c-a$ \\ \hline
\end{tabular}%
.
\end{equation}
\end{center}

On his page 3 Sallows uses $a=3,b=1,c=5$ to obtain the $Luoshu$ in (1).
However Sallows referenced neither Frierson, nor Andrews. See also Lachal%
\cite{Lachal}.

We tested $Lucas(3,1,5)$ finding numerical SVs $15,6.\,\allowbreak
928\,2..,3.\,\allowbreak 464\,1..$ $\allowbreak \allowbreak $which agree
with those of $M_{3}$, as expected.

\section{Other Magic Squares}

\subsection{Order $4$}

At order $4$ there are $880$ distinct magic squares of $1,2,..,16$ which
have been classified by the patterns of complementary number pairs within
the square into $12$ Groups by Dudeney\cite{Dudeney}. Counted in 1693 by Fr%
\'{e}nicle de \ Bessy, amongst them $48$\ associative and another $48$\ of
the pandiagonal variety defined soon. Since our 2013 study\cite{CRL} the 880
are now known to have $63$ different singular values (SV) clans, some of
which have just one non-zero EV\cite{LCR}.

\subsubsection{Pandiagonal Magic Squares (PMSs)}

Of the $880$, the $48$ in Dudeney\cite{Dudeney} Group I are pandiagonal.
These are characterised by having all parallel broken diagonals to the main
ones with the same RCD linesum as noted earlier for $AT_{3}$, but we note
that this is not the case for the present study of Frierson's associative
compound squares which we are not pandiagonal\cite{Weisstein}. There are
also $16$ ultramagic squares with both the associative and pandiagonal
features.

However this does not rule out other magic squares of orders $9$, $27$, $81$%
,... from being pandiagonal, some are known and one noted below, and others
we could construct by compounding.

\subsubsection{Ultramagic Squares}

These have both the associative and pandiagonal properties and begin at
order $5\cite{LAA}$.

\section{A low entropy order $4$ Most-Perfect Pandiagonal [MPPD] Magic Square%
}

These MPPDs are found at order $4$ and multiples of order $4$. Here we
consider one of this variety from the classic study of Dame Kathleen
Ollerenshaw and David Br\'{e}e\cite{Olly}\cite{OllyBree}, from their cover
but here using the classic elements $1,2,..n^{2}$ instead of zero-based:

\begin{center}
$MPPD_{4\alpha }=%
\begin{bmatrix}
1 & 15 & 4 & 14 \\ 
8 & 10 & 5 & 11 \\ 
13 & 3 & 16 & 2 \\ 
12 & 6 & 9 & 7%
\end{bmatrix}%
,$
\end{center}

This has the lowest entropy of the MSs in the $3$ pandiagonal clans Dudeney%
\cite{Dudeney} Groups $1,2,3$, the Alpha clan\cite{CRL}, with $\lambda
_{i}=34$, $\pm 8,0$, and $\sigma _{1}=$ $34.0,17.\,\allowbreak
889..,4.\,\allowbreak 472\,1..,0\allowbreak \allowbreak $, rank $3$.

Ollerenshaw \& Br\textit{\'{e}}e did not study any spectra, nor did Bellew,
but the former did reference Bellew.

\subsection{Compounding the lowest entropy $MPPD_{4\protect\alpha }$ for
comparison with our $n=3^{l}$ CMSs}

\begin{center}
\begin{tabular}{|l|l|l|l|l|}
\hline
$l$ & $1$ & $2$ & $3$ & $4$ \\ \hline
$n=4^{l}$ & $4$ & $16$ & $64$ & $256$ \\ \hline
$\sigma _{1}=\lambda _{1}$ & $34$ & $2056$ & $131,104$ & $8,388,736$ \\ 
\hline
$\sigma _{2}/\sqrt{5}$ & $2^{3}$ & $2^{9}$ & $2^{15}$ & $2^{21}$ \\ \hline
$\sigma _{3}/\sqrt{5}$ & $2$ & $2^{7}$ & $2^{13}$ & $2^{19}$ \\ \hline
$\sigma _{4}/\sqrt{5}$ &  & $2^{5}$ & $2^{11}$ & $2^{17}$ \\ \hline
$\sigma _{5}/\sqrt{5}$ &  & $2^{3}$ & $2^{9}$ & $2^{15}$ \\ \hline
$\sigma _{6}/\sqrt{5}$ &  &  & $2^{7}$ & $2^{13}$ \\ \hline
$\sigma _{7}/\sqrt{5}$ &  &  & $2^{5}$ & $2^{11}$ \\ \hline
$\sigma _{8}/\sqrt{5}$ &  &  &  & $2^{9}$ \\ \hline
$\sigma _{9}/\sqrt{5}$ &  &  &  & $2^{7}$ \\ \hline
$C\%$ & $37.2284..$ & $64.3023..$ & $75.9175..$ & $81.9199..$ \\ \hline
$H$ & $0.8702..$ & $0.98975..$ & $1.00156..$ & $1.00257..$ \\ \hline
$r$ & $3$ & $5$ & $7$ & $9$ \\ \hline
\end{tabular}

Table 8 - A lowest entropy order 4 magic square compounded.
\end{center}

In Table 8 the $\sigma _{2,3,4,5,..}$ increase by factors of $64$ across
columns from left to right as $l$ increases, and the $\sigma _{3,4,5,..}$\
decrease by factors of $1/4$ from their $\sigma _{2}$'s. The trend to an
asymptotic entropy mirrors that found in the main text for the lowest
entropy members of the Frierson CMSs.

\subsection{Higher orders $n=5,6,...$}

The populations of larger MSs continue to grow - see our colleague Walter
Trump's \ table \cite{Trump} which is regularly updated - so that the number
of distinct order $9$ MSs is astronomical, meaning that our Frierson-type
CMS are rare, but possibly close to the lowest entropy member?

\section{Numerical Compounding of Doubly Affine Matrices}

From c. 2004 Rogers and Cameron explored the use of Kronecker products of
MSs to generate larger ones of compound order - this was finally published
in 2017 with Loly\cite{RCL}[RCL]. RCL gave a general study of CMSs which
included these ancient pairs for arbitrary $m,n>2$ in terms of Kronecker
products, including a full account of their spectral properties which showed
that all CMSs are singular, a feature realized by them from earlier matrix
eigenvalue studies c. 2004. This was first reported at IWMS-2007, but not
included in the conference proceedings\cite{LAA} in 2009. RCL used the
entropy and compression measures from their 2013 CRL\cite{CRL} for magic and
Latin squares, including Frierson's first order 9 pair. Mixing orders $3$
and $4$ yields order $12$ compound squares with either order $4$ or $3$
subsquares, and is extendable to high orders.

RCL contains much useful background to compounding that need not be repeated
here as our focus is Frierson's different parameterized method, however one
result useful in the present context is that the rank of a CMS is the sum of
the ranks of its two components less $1$, e.g. the \ rank of a Frierson
order $9$ CMS is $3+3-1=5$, and for an order $27$ is then $5+3-1=7$, in
agreement with Table 5.

Note that our terminology for CMSs of $T_{A,D}$\ used here\ was changed in
RCL\cite{RCL} to\ $C_{A,D}$\ which are vertical reflections using the
sequence\ $0,1,2,..(n^{2}-1)$.

\section{Other tools for spectral calculations}

The authors have used Mathematica\cite{Mathematica} and Maple\cite{Maple}
software for both numerical and algebraic calculations as well as Numpy and
Sympy from Python\cite{Python}. In previous studies with Adam Rogers\cite%
{RCL}, also MATLAB, which has a magic square generator for one of each odd,
even and doubly-even orders. Earlier in "Online tools for calculating SVs"
we noted Keisan and Bluebit. We add that Wolfram Alpha\cite{Alpha} enables
access to some of Mathematica online, and well as via apps for iPhones and
iPads.\textbf{\ }

Also this article has been edited with a version of LaTeX in Scientific
Workplace\cite{SWP}, which also has a "Compute" section using Maple which
has been used recently to check some of the matrices herein.

\subsection{From the SVs to $C\%$ - a sample entropy and Compression
calculation for $M_{3}$}

The (default) numerical precision in SWP's\cite{SWP} "Evaluate Numerically"
is used here, first the total sigmas:

$\allowbreak 15+6.\,\allowbreak 928\,2+\allowbreak 3.\,\allowbreak
464\,1=25.\,\allowbreak 392$

then the contributions to the Shannon \ entropy, $H$, are calculated:

$-15/25.\,\allowbreak 392\times$ $\ln [15/25.\,\allowbreak 392]=\allowbreak
0.310\,95$

$-\allowbreak 6.\,\allowbreak 928\,2/25.\,\allowbreak 392\times\ln
[\allowbreak 6.\,\allowbreak 928\,2/25.\,\allowbreak 392]=\allowbreak
0.354\,39$

$-\allowbreak 3.\,\allowbreak 464\,1/25.\,\allowbreak 392\times\ln
[\allowbreak 3.\,\allowbreak 464\,1/25.\,\allowbreak 392]=\allowbreak
0.271\,76$

For a total: $H=0.310\,95+0.354\,39+0.271\,76=\allowbreak 0.937\,1$, and
finally the \% Compression follows:

$C=(1-\allowbreak 0.937\,1/\ln [\allowbreak 3])\times100=\allowbreak
14.\,\allowbreak 701$, which both agree with our CRL\cite{CRL} calculations.

\subsection{A Cautionary note}

Since some computer software, e.g \ MATLAB\cite{MatLab} and "R"\cite{R},
offer just a single magic square for each order one must be careful to not
draw strong conclusions from their single MSs as to the properties of others
of the same order in view of the great diversity already apparent at order $%
4 $. Clearly our Frierson-type associative CMSs are going to be just a
(small) fraction of the enormous number of order $9$ magic squares, but
perhaps of low entropy.

\section{Browne's CMS27}

An order $27$\ CMS by Browne\cite{Browne}, $B27$, with a commentary by Paul
Carus, was shown in chapter VI of MSC1\cite{Andrews}, Fig. 273 (Fig.256 of
MSC2\cite{Andrews}), but is not easy to read, in\ part because alternate
cells are shaded.

$B27$ may be a variant of our \#$15$, "O" in Table 7 [$%
v=27,y=1;s=81,t=3;p=243,q=9$]\textbf{. }

For compactness and accuracy we divide $F_{27}$ into $9$-by-$9$ order $9$
subsquares and those similarly to order $3$ subsquares beginning with:

$B3=%
\begin{bmatrix}
28 & 57 & 2 \\ 
3 & 29 & 55 \\ 
56 & 1 & 30%
\end{bmatrix}%
$, for which $v=27,y=1$,

and which has Knight path's $1\rightarrow 2$ and $2\rightarrow3$, then a
move up by $25$\ to begin a down diagonal $28\rightarrow29\rightarrow30$,
and two more Knight path's $55\rightarrow56$ and $56\rightarrow57$.

Then using \ to construct the bottom middle order 9 subsquare:

$B9=%
\begin{bmatrix}
B3+81E3 &  & B3+168E3 &  & B3+3E3 \\ 
B3+6E3 &  & B3+84E3 &  & B3+162E3 \\ 
B3+165E3 &  & B3 &  & B3+87E3%
\end{bmatrix}%
$,

and finally:

$B27=%
\begin{bmatrix}
B9+243E9 &  & B9+504E9 &  & B9+9E9 \\ 
B9+18E9 &  & B9+252E9 &  & B9+486E9 \\ 
B9+495E9 &  & B9 &  & B9+261E9%
\end{bmatrix}%
,$

but not given explicitly as it takes a whole page - see the clarity issue in
Browne's\cite{Browne} example in MSC2\cite{Andrews}, page 150, which is
clearer in Swetz\cite{Swetz}, page 136, with a duplicate 606 in row 7,
column 15 which should be 506.

\bibliographystyle{aaai-named}
\bibliography{acompat,testAPA,texbook1}

\begin{thebibliography}{99}
\bibitem{Andrews} [MSC1,2] Andrews, W.S., \textit{Magic Squares and Cubes},
[MSC1 1908], The Open Court Publishing Corporation 1908 with particularly
notable chapters by L.S. Frierson and C.A. Browne. Introduction, vii-viii,
by Dr. Paul Carus, editor of \textit{The Monist} (see also Browne\cite%
{Browne} and Frierson\cite{Frierson}). For full text -\TEXTsymbol{>}
https://archive.org/?z. See also Andrews, W.S., \textit{Magic Squares and
Cubes}, [MSC2 1917], The Open Court Publishing Corporation 1917, Reprinted
by Cosimo 2004.

\bibitem{Bellew} Bellew, J., 1997, \textit{Counting the Number of Compound
and Nasik Magic Squares}, \textit{Mathematics Today}, \textbf{33}(4),
111--118.

\bibitem{Bergholt} Bergholt, E., \textit{The Magic Square of Sixteen Cells.
A New and Completely General Formula}, Nature \textbf{83}(1910), No. 2117,
pp 368-9.

\bibitem{Browne} Browne, C.A., in \cite{Andrews} see Ch. VI, \textit{Magics
and Pythagorean Numbers,} see p. 151 for n=27 - see footnote pp. 158-162 by
Paul Carus.

\bibitem{bluebit} Bluebit Software, \textit{Linear algebra math tools for
.NET, C\#, VB, C++} , https://pressaboutus.com/bluebit.gr

\bibitem{CRL} [CRL] Cameron, I.D. Rogers, A. and Loly, P.D. \textit{%
Signatura of magic and Latin integer squares: isentropic clans and indexing}%
, Discuss. Math. Probab. Stat., \textbf{33} (2013), pp. 121-149,
http://www.digora.pl/, from 2012 Bedlewo conference. Download paper:from
http://www.discuss.wmie.uz.zgora.pl/ps. PowerPoint talk from 2012 conference
paper published in Discussiones Mathematicae Probability and Statistics,
33(1-2) (2013) 121-149.-\TEXTsymbol{>} http://home.cc.umanitoba.ca/\symbol{%
126}loly/Bedlewo.txtvideo. - \ PowerPoint: -\TEXTsymbol{>}
http://home.cc.umanitoba.ca/\symbol{126}loly/Signatura.pdf; data advertised
but no link - see \cite{CRLdata} next:

\bibitem{CRLdata} [CRLdata] Cameron, I.D. Rogers, A. and Loly, P.D.,
\textquotedblleft Data Appendix for Signatura\textquotedblright ,
DMPS2013http://home.cc.umanitoba.ca/\symbol{126}loly/Bedlewo.txt

\bibitem{Cameron} Cameron, I.D., provisional title: \textit{Four Parameter
Study of Particular Fr\'{e}nicle de Bessy's Magic Squares}, in preparation.

\bibitem{Cammann} Cammann, Schuyler, 1960, \textit{The evolution of magic
squares in China}, Journal of the American Oriental Society, \textbf{80}
(2), 116--124; (1962) \textit{Old Chinese Magic Squares}, Sinologica, 
\textbf{7}, 14-53 ; Schuyler Cammann (1957). \textit{Magic square}, Encyclop%
\ae dia Britannica, 14th edition, vol. XIV, page 625 TBC. Revised 14th
Edition, vol. 14, 573-5. [revised 14th edition 1933-1973; Schuyler Cammann
(1961), \textit{The magic square of three in old Chinese philosophy and
religion}, History of Religions, \textbf{1} (1), 37--80.

\bibitem{Chan} Chan, W. and Loly, P. D., [CL] 2002, \textit{Iterative
Compounding of Square Matrices to Generate Large-Order Magic Squares},
Mathematics Today, \textbf{38}(4), 113-118, (The Institute of Mathematics
and its Applications, Southend-on-Sea, UK).

\bibitem{Derksen} Derksen, Harm, Christian Eggermont, and Arno van den Essen,%
\textit{\ Multimagic Squares},https://arxiv.org/abs/math/0504083

\bibitem{Descombes} Descombes, R., 2000, \textit{Les Carr\'{e}s Magiques:
Histoire, th\'{e}orie et technique du carr\'{e} magique, de l'Antiquit\'{e}
aux recherches actuelles,} Vuibert (Paris).

\bibitem{Drury} Drury, S., \textit{There are no magic squares of rank 2},
personal communication c.2007.

\bibitem{Dudeney} Dudeney, H.E. in Encyclopedia Britannica, 1867, article on 
\textit{magic squares} [also 1929 in 14th edition], and (now with) R. C.
Bo(se) in 15th edition 1974.

\bibitem{Frierson} Frierson, L.S., \textit{A Mathematical Study of Magic
Squares: A New Analysis}, The Monist, \textbf{XVII} (1907), 272-293, (in
Criticism and Discussion, signed L.S.Frierson, Frierson, LA.). This is
reproduced in Andrews\cite{Andrews} as Chapter 5, pp. 129-137.

\bibitem{Heinz} Heinz, H. D. and Hendricks, J. R. 2000, \textit{Magic
Squares Lexicon: Illustrated},
http://magic-squares.net/Downloads/HendricksBooks/Lexicon-v2.pdf

\bibitem{HH} Heinz, H.D., \textit{Order 4 magic-squares} -\TEXTsymbol{>}
magic-squares.net/order4list.htm

\bibitem{Hendricks} Hendricks, J. R., 1992, \textit{The Magic Square Course}%
, 2nd edition, (privately printed). Many of Hendricks results have been
included in Pickover \cite{Pickover}.

\bibitem{HJ} Horn, R.A. and Johnson, C.R., \textit{Matrix Analysis}, Second
Edition, Cambridge University Press, 2013.

\bibitem{Keisan} Keisan - \textit{Online Singular Value Decomposition
Calculator,} https://keisan.casio.com/exec/system/15076953160460 \ - enter
matrix - see SVs in middle row "Wj"

\bibitem{Kirkland} Kirkland, S.J. and Neumann, M., \textit{Group Inverses of
M-matrices associated with non-negative matrices having few eigenvalues},
Linear Algebra and Its Applications, 220 (1995) 181-213.

\bibitem{Lachal} Lachal, Aim\'{e}, \textit{Carr\'{e}s Magiques a l'ordre 3}
par Lucas, \'{E}douard - link: {\small math.univ-lyon1.fr/\symbol{126}%
alachal/diaporamas/diaporama\_carres\_magiques\_ordre3.pdf \& }%
math.univ-lyon1.fr/\symbol{126}alachal/exposes/les\_carres\_magiques1.pdf,
ibid /\symbol{126}alachal/exposes/carres\_magiques\_diaporama.pdf

\bibitem{Lam} Lam Lay Yong (1977), \textit{A Critical Study of the Yang Hui
Suan Fa : A Thirteenth-century Chinese Mathematical Treatise,} Singapore
University Press. [\textquotedblleft This book is an amended version of a
thesis submitted for the degree of Ph.D. at the University of
Singapore.\textquotedblright ]

\bibitem{LiYen} Li Yen (1954). \textit{Chung Suan Shih Lun Ts'ung = A
Discussion on the History of Chinese Mathematics}, 5 vols. Science
Publishing Company, Peking. Li Yen (1990). Chung Suan Shih Lun Ts'ung = A
Discussion on the History of Chinese Mathematics, 1282 pp. \textit{Shang-hai
shu tien}, Shanghai. [OCLC: 70359520 at Univ. Alberta. Apparently reprint
edition of Li Yen (1934/1936).

\bibitem{LAA} [LAA] Loly, P., Cameron, I., Trump, W. and Schindel, D., 
\textit{Magic square spectra}, Linear Algebra and its Applications, \textbf{%
430} (2009) 2659-2680.

\bibitem{Loly2009} [LC] Loly, P. (with I.D.Cameron), \textit{Eigenvalues of
an Algebraic Family of Compound Magic Squares of Order} $n=3^{l},l=1,2,3,..$%
, \textit{and Construction and Enumeration of their Fundamental Forms},
PowerPoint presentation for lead keynote talk by Loly at Canadian
Mathematical Society 2009, Windsor, Ontario.

\bibitem{LCR} [LCR] Loly, P.D., Cameron, I.D. and Rogers, A., \textit{Powers
of doubly-affine integer square matrices with one non-zero eigenvalue},
arXiv:1712.03393[math.HO]

\bibitem{Lucas} Lucas, \'{E}douard, \textit{Carr\'{e}s Magiques a l'ordre 3}%
, Recreations Mathematiques \textbf{IV} (1894):225 \cite{Lachal}

\bibitem{Mandelbrot} Mandelbrot, B. B., \textit{The Fractal Geometry of
Nature},\textit{\ }Macmillan, 1983 ISBN 978-0-7167-1186-5

\bibitem{Maple} Char, B.W. and 5 others, \textit{Maple V Library Reference
Manual}, Springer-Verlag, 1991.

\bibitem{Mathematica} \textit{Mathematica: Modern Technical Computing},
www.wolfram.com/mathematica

\bibitem{MatLab} \textit{MATLAB }%
https://www.mathworks.com/products/matlab.html

\bibitem{NDS} [NDS] Newton, P.K. and DeSalvo, S.S. \textit{The Shannon
entropy of Sudoku matrices}, Proc. R. Soc. \textbf{466} (2010), 1957-1975.

\bibitem{OllyBree} Ollerenshaw, K. and Br\textit{\'{e}}e, D. S., 1998, 
\textit{Most-perfect pandiagonal magic squares: their construction and
enumeration, }The Institute of Mathematics and its Applications,
Southend-on-Sea, UK.

\bibitem{Olly} Ollerenshaw, K., 2006, \textit{Constructing pandiagonal magic
squares of arbitrarily large size}, Mathematics Today, \textbf{42,} Parts 1
and 2, Feb. 23-29; Part 3, Apr. 66-69, The Institute of Mathematics and its
Applications, Southend-on-Sea, UK.

\bibitem{Pickover} Pickover, C., 2002, \textit{The Zen of Magic Squares,
Circles, and Stars - An Exhibition of Surprising Structures across
Dimensions,} Princeton University Press, Princeton, New Jersey.

\bibitem{Python} python.org for Numpy and Sympy.

\bibitem{R} \textit{The R Project for Statistical Computing}
https://www.r-project.org

\bibitem{RCL} [RCL] Rogers, A., Cameron, I.D. and Loly, P.D. [RCL], \textit{%
Compounding Doubly Affine Matrices,} arXiv:1711.11084 [math.CO] 2017, 
\textit{arXiv}:1711.11084[math.CO]

\bibitem{Sallows} Sallows, Lee C.F., \textit{Geometric Magic Squares - A
Challenging New Twist Using Colored Shapes Instead of Numbers}, Dover
Publications, Inc. New York.

\bibitem{Schubert} Hermann Schubert (1898), \textit{Mathematical Essays and
Recreations}, translated from German into English by Thomas J. McCormack, 
\textit{Open Court Publishing Corporation}, Chicago - {\tiny %
http://ia350605.us.archive.org/1/items/mathematicalessa00schuiala/mathematical?essa00schuiala.pdf%
}{\scriptsize \texttt{\ }}

\bibitem{Shreveport} Frierson, Louisiana - Wikipedia, Frierson, Louisiana -
Wikipedia: https://en.wikipedia.org/wiki/Frierson,\_Louisiana

\bibitem{Sesiano} Sesiano, J., 2004, \textit{Les Carr\'{e}s Magiques dans
les Pays Islamic}, Presses Polytechniques et Universitaires Romandes;
Sesiano, Jacques, 1996, \textit{Trait\'{e} M\'{e}di\'{e}val sur les Carr\'{e}%
s Magiques}, Presses Polytechniques et Universitaires Romandes, Sesiano,
Jacques, 1998, \textit{Le Trait\'{e} d' Ab\={u}'l-Waf\={a}' sur les Carr\'{e}%
s Magiques}, f\"{u}r Geschichte der Arabisch-Islamischen wissenschaften,
Band 12, 121-200-244.

\bibitem{Sloane384} Sloane, N.J.A., \textit{Hexagonal numbers}: $n(2n-1)$.
The On-Line Encyclopedia of Integer Sequences. A000384
http://www.research.att.com/\symbol{126}njas/sequences/A000384

\bibitem{Sloane680} Sloane, N.J.A.,$\ \left( 2n\right) !/2\symbol{94}n$ The
On-Line Encyclopedia of Integer Sequences. A000680 {\small %
http://www.research.att.com/\symbol{126}njas/sequences/A000680}

\bibitem{Sloane1147} Sloane, N.J.A.,$\ a(n)=(2n-1)!!\ $\textit{Double
factorial of odd numbers}, The On-Line Encyclopedia of Integer Sequences.
A0001147 http://www.research.att.com/\symbol{126}njas/sequences/A0001147

\bibitem{Sloane523386} Sloane, N.J.A., The On-Line Encyclopedia of Integer
Sequences. A0523386 http://www.research.att.com/\symbol{126}%
njas/sequences/A0523386

\bibitem{Swetz5} Swetz, F.J., 1994, \textit{From Five Fingers to Infinity -
A Journey through the History of Mathematics, }edited by Frank J. Swetz,
Open Court\textit{.}

\bibitem{Swetz} Swetz, F. J., 2002, \textit{Legacy of the Luoshu - The 4000
Year Search for the Meaning of the Magic Square of Order Three,} Chicago:
Open Court; see also Swetz (1977), \textit{The `piling up of squares' in
ancient China}, The Mathematics Teacher, \textbf{70}, 72--79; Swetz (1978)
and \textbf{71} (1), 50--56 (January 1978).

\bibitem{SWP} [SWP] Hardy, D.W. and Walker, C.L., \textit{Doing Mathematics
with Scientific WorkPlace and Scientific Notebook}, v. 5.5, 2005.

\bibitem{Thompson} Thompson, W.H., \textit{On Magic Squares}, The Quarterly
Journal of Pure and Applied Mathematics (1869) 186-202.

\bibitem{Trigg} Trigg, C. W., 1980, \textit{A Family of Ninth Order Magic
Squares}, Mathematics Magazine, \textbf{53}, 110-1.

\bibitem{Trump} Trump, W.,\textit{\ Notes on Magic Squares and Cubes}, 2003,
www.trump.de/magic-squares/howmany.html

\bibitem{Violle} Violle, B., 1837, \textit{Trait\'{e} complet des carr\'{e}s
magiques} (1837---1838) https://www.amazon.ca/Magiques-Compos\'{e}%
s-Bordures-Compartimens-D\'{e}tach\'{e}es/dp/0266470076

\bibitem{Weisstein} Weisstein, E., \textit{Associative Magic Square} on
Wolfram MathWorld: https://mathworld.wolfram.com/AssociativeMagicSquare.html

\bibitem{WeissteinPan} Weisstein, E., \textit{Pandiagonal Magic Square} on
Wolfram MathWorld: https://mathworld.wolfram.com/PanmagicSquare.html

\bibitem{White} White, S. Harry: budshaw.ca/index.html

\bibitem{Alpha} Wolfram Alpha: https://www.wolframalpha.com/
\end{thebibliography}

\end{document}